\documentclass[10pt]{article}
\usepackage{latexsym}

\newtheorem{veta}{Theorem}[section]
\newtheorem{dusl}[veta]{Corollary}
\newtheorem{lema}[veta]{Lemma}
\newtheorem{defi}[veta]{Definition}

\newtheorem{prop}[veta]{Proposition}

\newcommand{\duk}{\noindent {\bf Proof. }}
\newcommand{\kduk}{\hfill $\Box$\bigskip}
\newcommand{\N}{\mathbf{N}}
\newcommand{\Z}{\mathbf{Z}}

\newcommand{\M}{{\cal S}}

\author{Martin Klazar\thanks{Institute for Theoretical Computer Science and Department of Applied Mathematics, 
Faculty of Mathematics and Physics of Charles University, Malostransk\'e n\'am\v est\'\i\ 25, 118 00 Praha, Czech Republic. 
ITI is supported by the project 1M0021620808 of the Czech Ministry of Education.
Email: {\tt klazar@kam.mff.cuni.cz}}}
\title{On identities concerning the numbers of crossings and nestings of two edges in matchings}
\date{\today}

\begin{document}

\maketitle
\begin{abstract}
Let $M,N$ be two matchings on $[2n]$ (possibly $M=N$) and for an integer $l\ge 0$ let ${\cal T}(M,l)$ be the set 
of those matchings on $[2n+2l]$ which can be obtained from $M$ by successively adding $l$ times in all ways 
the first edge, 
and similarly for ${\cal T}(N,l)$. Let $s,t\in\{cr,ne\}$ where $cr$ is the statistic of the number of crossings 
(in a matching) and 
$ne$ is the statistic of the number of nestings (possibly $s=t$). We prove that if the statistics $s$ and $t$ 
coincide on the sets of matchings ${\cal T}(M,l)$ and ${\cal T}(N,l)$ for $l=0,1$, they must coincide on these sets
for every $l\ge 0$; similar identities hold for the joint statistic of $cr$ and $ne$. These results are instances of 
a general identity in which crossings and nestings are weighted by elements from an abelian group.
\end{abstract}

\section{Introduction and formulation of the main result}

In this article we investigate distributions of the numbers of crossings and nestings of two edges in matchings. For example, 
it is known that for each $k$ and $n$ there are as many matchings $M$ on $\{1,2,\dots,2n\}$ with $k$ crossings as 
those with $k$ nestings. All matchings form an infinite tree ${\cal T}$ rooted in the empty matching $\emptyset$, 
in which the children 
of $M$ are the matchings obtained from $M$ by adding to $M$ in all possible ways new first edge. The problem we address is this: 
Given two (not necessarily distinct) matchings $M$ and $N$ on $\{1,2,\dots,2n\}$, when is it the case that the numbers of crossings
(or nestings, or crossings versus nestings) have the same distributions on the levels of the two subtrees of ${\cal T}$ rooted in $M$
and $N$. Our main result is Theorem~\ref{generalthm} that determines when this happens, in fact in a more general setting. 
Before formulating it we give definitions and fix notation.

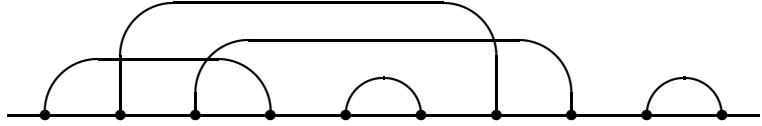
\begin{figure}\label{parovani}
\unitlength1mm
\begin{picture}(30,20)(-2,-3)
\put(20,0){\line(1,0){100}}
\multiput(25,0)(10,0){10}{\circle*{1.5}}
\thicklines
\put(40,0){\oval(30,15)[t]}
\put(60,0){\oval(50,30)[t]}
\put(70,0){\oval(50,20)[t]}
\put(70,0){\oval(10,10)[t]}
\put(110,0){\oval(10,10)[t]}
\end{picture}
\caption{Matching with 3 crossings and 2 nestings.}
\end{figure}

We denote the set $\{1,2,3,\dots\}$ by $\N$, the set $\N\cup\{0\}$ by $\N_0$, and (for $n\in\N$) the set $\{1,2,\dots,n\}$ by $[n]$. 
The cardinality of a set $A$ is denoted $|A|$. By a {\em multiset} we understand a ``set'' in which repetitions of elements are 
allowed. This can be modeled by a pair $H=(X,m)$ where $X$ is a set, the {\em groundset} of the multiset $H$, and the mapping 
$m:\;X\to\N$ determines the multiplicities of the elements in $H$. However, we will not need this formalism and will record 
multiplicities by repetitions. A {\em matching} $M$ on $[2n]$ is a set partition of $[2n]$ in $n$ two-element 
{\em blocks} which we also call  
{\em edges}. The set of all matchings on $[2n]$ is denoted ${\cal M}(n)$; we define ${\cal M}(0)=\{\emptyset\}$. Two distinct blocks 
$A$ and $B$ of $M$ form a {\em crossing} (they {\em cross}) if $\min A<\min B<\max A<\max B$ or $\min B<\min A<\max B<\max A$. 
Similarly, they form a {\em nesting} (they are {\em nested}) if $\min A<\min B<\max B<\max A$ or 
$\min B<\min A<\max A<\max B$. We draw a diagram of $M$ in which we put the elements $1,2,\dots,2n$ as points on a line, from left 
to right, and connect by a semicircular arc lying above the line the two points of each block. 
For two crossing blocks the corresponding 
arcs intersect and for two nested blocks one of the arcs covers the other, see Figure~\ref{parovani}. 
By $cr(M)$, respectively $ne(M)$, we denote the number of crossings, respectively nestings, in $M$. The $n$ edges of 
$M\in{\cal M}(n)$ are naturally ordered by their first elements. The first edge of $M$ is $\{1,x\}$ and the last edge is 
the one whose first vertex is the last one among the $n$ first vertices.

We investigate distribution of the numbers $cr(M)$ and $ne(M)$ on ${\cal M}(n)$ and on the subsets of ${\cal M}(n)$ defined
by prescribing the matching formed by the last $k$ edges of $M$. The total number of matchings in ${\cal M}(n)$ is 
$$
|{\cal M}(n)|=(2n-1)!!=1\cdot 3\cdot 5\cdot\dots\cdot(2n-1). 
$$
It is known that the number of matchings on $[2n]$ with no crossing equals the number of matchings with no nesting and that it is 
the $n$-th Catalan number, see Stanley \cite[Problems 6.19o and 6.19ww]{stan}:
$$
|\{M\in{\cal M}(n):\;cr(M)=0\}|=|\{M\in{\cal M}(n):\;ne(M)=0\}|=\frac{1}{n+1}{2n\choose n}.
$$
The more general result that for each $k$ and $n$
$$
|\{M\in{\cal M}(n):\;cr(M)=k\}|=|\{M\in{\cal M}(n):\;ne(M)=k\}|
$$
was derived by M. de Sainte-Catherine in \cite{desa}. Even more is true because the joint statistic is symmetric:
$$
|\{M\in{\cal M}(n):\;cr(M)=k,ne(M)=l\}|=|\{M\in{\cal M}(n):\;cr(M)=l,ne(M)=k\}|
$$
for every $k,l\in\N_0$ and $n\in\N$. A simple proof for this symmetry can be given by adapting the Touchard-Riordan method 
(\cite{touc}, \cite{rior}) that encodes matchings and their numbers of crossings by weighted Dyck paths, see Klazar and 
Noy \cite{klaz_noy}. Here we put these results in a more general framework.

\begin{figure}
\unitlength1mm
\begin{picture}(60,60)(0,0)
\put(69.3,0){$\emptyset$}
\thicklines
\put(70,5){\line(0,1){6}}
\put(70,13){\oval(8,6)[t]}
\put(62,15){\line(-3,1){21.5}}
\put(70,18){\line(0,1){6}}
\put(78,15){\line(3,1){21.5}}
\put(34,26){\oval(8,9)[t]}\put(46,26){\oval(8,9)[t]}
\put(68,26){\oval(8,9)[t]}\put(72,26){\oval(8,9)[t]}
\put(100,26){\oval(8,6)[t]}\put(100,26){\oval(16,12)[t]}
\put(40.5,35){\line(0,1){8}}
\put(38.5,35){\line(-2,3){5}}
\put(42.5,35){\line(2,3){5}}
\put(35.5,35){\line(-4,3){7}}
\put(45.5,35){\line(4,3){7}}
\put(100.5,35){\line(0,1){8}}
\put(98.5,35){\line(-2,3){5}}
\put(102.5,35){\line(2,3){5}}
\put(95.5,35){\line(-4,3){7}}
\put(105.5,35){\line(4,3){7}}
\put(70.5,35){\line(0,1){8}}
\put(68.5,35){\line(-2,3){5}}
\put(72.5,35){\line(2,3){5}}
\put(65.5,35){\line(-4,3){7}}
\put(75.5,35){\line(4,3){7}}
\put(40,46){$\vdots$}
\put(70,46){$\vdots$}
\put(100,46){$\vdots$}
\end{picture}
\caption{Tree of matchings ${\cal T}$.}
\label{st}
\end{figure}
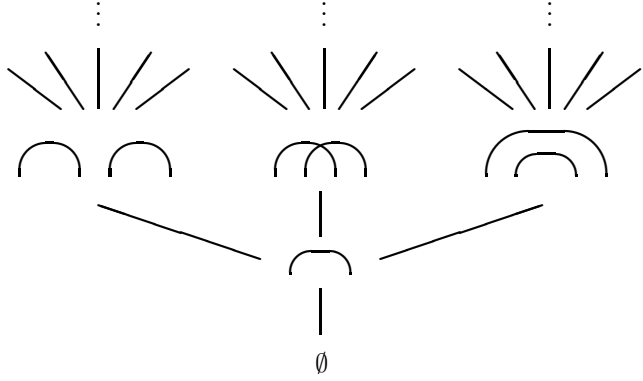

By the {\em tree of matchings} ${\cal T}=({\cal M},E,r)$ we understand the infinite rooted tree with 
the vertex set 
$$
{\cal M}=\bigcup_{n=0}^{\infty}{\cal M}(n),
$$
which is rooted in the empty matching $r=\emptyset$ and in which 
directed edges in $E$ are the pairs $(M,N)$ such that $M\in{\cal M}(n)$, $N\in{\cal M}(n+1)$, and 
$N$ arises from $M$ by adding a new first edge, that is, we relabel the vertices of $M$ as 
$\{2,3,\dots,2n+2\}\backslash\{x\}$ for some  $x\in\{2,3,\dots,2n+2\}$ and add to $M$ the  
block $\{1,x\}$, see Figure~\ref{st}. 
Each vertex $N\in{\cal M}(n)$ has $2n+1$ children and, if $n>0$, is a 
child of a unique vertex $M\in{\cal M}(n-1)$. A {\em level} in a rooted tree is the set of vertices with the 
same distance from the root. In ${\cal T}$ the levels are the sets ${\cal M}(n)$. The 
{\em subtree} ${\cal T}(M)$ of ${\cal T}$ rooted in $M\in {\cal M}(n)$ is the rooted subtree on the vertex 
set ${\cal N}\subset {\cal M}$ consisting of $M$ and all its descendants, that is, ${\cal N}$ contains $M$ and all 
matchings obtained from $M$ by successively adding new first edge. In other words, ${\cal T}(M)$ consists of all 
$N\in{\cal M}$ in which the last $n$ edges form a matching (order-isomorphic to) $M$. Clearly, ${\cal T}(\emptyset)={\cal T}$. 
We denote the $l$-th level of ${\cal T}(M)$ by ${\cal T}(M,l)$. For $M\in{\cal M}(n)$ we have ${\cal T}(M,0)=\{M\}$ and 
${\cal T}(M,1)$ is the set of children of $M$ in ${\cal T}$. Also, $|{\cal T}(M,l)|=(2n+1)(2n+3)\dots(2n+2l-1)$.  

Besides the statistics $cr(M)\in\N_0$ and $ne(M)\in\N_0$ on ${\cal M}$ we consider the joint statistics 
$cn(M)=(cr(M),ne(M))\in\N_0^2$ and $nc(M)=(ne(M),cr(M))\in\N_0^2$.
Two statistics $s,u$ on two subsets ${\cal N}_1,{\cal N}_2\subset{\cal M}$ {\em coincide (have the same distribution)} if 
$s({\cal N}_1)=u({\cal N}_2)$ as multisets, that is, if for every element $e$ we have 
$$
|\{M\in{\cal N}_1:\;s(M)=e\}|=|\{M\in{\cal N}_2:\;u(M)=e\}|.
$$

\bigskip\noindent
{\bf Notational convention.} If $f:\;X\to Y$ is a mapping and $Z\subset Y$, the symbol $f(Z)$ usually denotes the image
$Im(f|Z)=\{f(z):\;z\in Z\}$. In this article we use $f(Z)$ to denote the multiset whose ground set is $Im(F|Z)$ and in which 
each element $y=f(z)$, $z\in Z$, appears with the multiplicity $|f^{-1}(y)\cap Z|$. So in our $f(Z)$ each element $y$ has the proper 
multiplicity in which it is attained as a value of $f$ on $Z$.

\bigskip
Let $A=(A,+)$ be an abelian group and $\alpha,\beta\in A$ be its two elements. The most general statistic on matchings that we 
consider is $s_{\alpha,\beta}:\;{\cal M}\to A$ given by
$$
s_{\alpha,\beta}(M)=cr(M)\alpha+ne(M)\beta.
$$
Our main result is the next theorem. 

\begin{veta}\label{generalthm}
Let $M,N\in{\cal M}(n)$ be two (not necessarily distinct) matchings and, for $\alpha,\beta\in A$, $s_{\alpha,\beta}$ be the 
above statistic.
\begin{enumerate}
\item If $s_{\alpha,\beta}({\cal T}(M,l))=s_{\alpha,\beta}({\cal T}(N,l))$ for $l=0,1$  
then $s_{\alpha,\beta}({\cal T}(M,l))=s_{\alpha,\beta}({\cal T}(N,l))$ for all $l\ge 0$.
\item If $s_{\alpha,\beta}({\cal T}(M,l))=s_{\beta,\alpha}({\cal T}(N,l))$ for $l=0,1$
then $s_{\alpha,\beta}({\cal T}(M,l))=s_{\beta,\alpha}({\cal T}(N,l))$ for all $l\ge 0$.
\end{enumerate}
In words, for the statistic $s_{\alpha,\beta}$ to coincide level by level on the subtrees ${\cal T}(M)$ and 
${\cal T}(N)$ it suffices if it coincides on the first two levels, and similarly for the pair of statistics 
$s_{\alpha,\beta},s_{\beta,\alpha}$.
\end{veta}

\noindent
Specializing, we obtain identities for the statistics $cr,ne,cn$, and $nc$.

\begin{veta}\label{mainthm}
Let $M,N\in{\cal M}(n)$ be two (not necessarily distinct) matchings and $s,t\in\{cr,ne\}$, 
$u,v\in\{cn,nc\}$ be statistics on matchings (we allow $s=t$ and $u=v$).
\begin{enumerate}
\item If $s({\cal T}(M,l))=t({\cal T}(N,l))$ for $l=0,1$ then $s({\cal T}(M,l))=t({\cal T}(N,l))$ for all $l\ge 0$.
\item If $u({\cal T}(M,l))=v({\cal T}(N,l))$ for $l=0,1$ then $u({\cal T}(M,l))=v({\cal T}(N,l))$ for all $l\ge 0$.
\end{enumerate}
\end{veta}
\duk
1. Let $A=(\Z,+)$. Setting $\alpha=1,\beta=0$ and $\alpha=0,\beta=1$ and using 1 and 2 of Theorem~\ref{generalthm}, we obtain 
the identities for $cr$ and $ne$. 

2. Let $A=(\Z^2,+)$. Setting $\alpha=(1,0),\beta=(0,1)$ and $\alpha=(0,1),\beta=(1,0)$ and using 1 and 2 of 
Theorem~\ref{generalthm}, we obtain the identities for $cn$ and $nc$.
\kduk

\noindent
We illustrate the last theorem by four examples. We mentioned the first two already, it is the result of de Sainte-Catherine and 
the symmetry $cn=nc$.

\begin{dusl}\label{desainte}
For every $k\in\N_0$ and $n\in\N$ there are as many matchings on $[2n]$ with $k$ crossings 
as those with $k$ nestings.
\end{dusl}
\duk
Set $M=N=\emptyset$ and $s=cr,t=ne$. The assumption of the theorem is satisfied
because $cr(\emptyset)=ne(\emptyset)=0$ and $cr({\cal M}(1))=ne({\cal M}(1))=\{0\}$. 
\kduk

\begin{dusl}\label{symmetric}
For every $k,l\in\N_0$ and $n\in\N$ there are as many matchings on $[2n]$ with $k$ crossings 
and $l$ nestings, as those with $l$ crossings and $k$ nestings; the joint statistic is symmetric.
\end{dusl}
\duk
Set $M=N=\emptyset$ and $s=cn,t=nc$. The assumption of the theorem is satisfied
because $cn(\emptyset)=nc(\emptyset)=(0,0)$ and $cn({\cal M}(1))=nc({\cal M}(1))=\{(0,0)\}$. 
\kduk

\begin{dusl}\label{crequalne}
For every $k\in\N_0$ and $n\in\N$ there are as many matchings on $[2n]$ which have $k$ crossings
and have the last two edges nested, as those  which have $k$ nestings and have the last two edges separated (neither crossing 
nor nested).
\end{dusl}
\duk
Set $M=\{\{1,4\},\{2,3\}\}$, $N=\{\{1,2\},\{3,4\}\}$, $s=cr$, and $t=ne$. The assumption of the theorem is satisfied 
because $cr(M)=ne(N)=0$ and the values of $cr$ on the five children of $M$ are $0,0,1,1,2$, which coincides 
with the values of $ne$ on the five children of $N$. 
\kduk
\begin{dusl}\label{crequalcr}
Let $M=\{\{1,2\},\{3,5\},\{4,6\}\}$ and $N=\{\{1,3\},\{2,4\},\{5,6\}\}$. For every 
$k,n\in\N$ there are as many matchings on $[2n]$ with $k$ crossings in which the last three edges form 
a matching order-isomorphic to $M$, as those in which the last three edges form a matching order-isomorphic to $N$. 
\end{dusl}
\duk
Set the matchings $M,N$ as given and $s=t=cr$. Then $cr(M)=cr(N)=1$ and 
$cr({\cal T}(M,1))=cr({\cal T}(N,1))=\{1,1,1,2,2,2,3\}$. 
\kduk

We call two matchings $M,N\in{\cal M}(n)$ {\em crossing-similar} and write $M\sim_{cr}N$ if 
$cr({\cal T}(M,l))=cr({\cal T}(N,l))$ for all $l\ge 0$. Similarly we define the {\em nesting-similarity} $\sim_{ne}$. 
These two relations are equivalences and partition ${\cal M}(n)$ in equivalence classes. We use Theorem~\ref{mainthm} to 
characterize these classes and to count them. In Theorems~\ref{cros_equi} and \ref{nest_equi} we prove that the numbers of 
classes in ${\cal M}(n)/\!\sim_{cr}$ and ${\cal M}(n)/\!\sim_{ne}$ are, respectively, 
$$
2^{n-2}\bigg({n\choose 2}+2\bigg)\ \mbox{ and }\ 
2\cdot 4^{n-1}-\frac{3n-1}{2n+2}{2n\choose n}.
$$
These two numbers differ, the latter is roughly a square of the former. On the first level of 
description of the enumerative complexity of crossings and nestings, that of the numbers $cr(M)$ and $ne(M)$, symmetry reigns 
as shown in Corollaries~\ref{desainte} and  \ref{symmetric}. On the next level of description, that of 
the similarity classes, symmetry is broken because $|{\cal M}(n)/\!\sim_{ne}\!|$ is much bigger than 
$|{\cal M}(n)/\!\sim_{cr}\!|$. From this point of view nestings are definitely more complicated than crossings, see also 
Theorem~\ref{modtwo}. 

We prove Theorem~\ref{generalthm} in Section 2. The method we employ is induction on the number of edges. In Section 3
we prove Theorems~\ref{cros_equi} and \ref{nest_equi} enumerating the crossing-similarity and nesting-similarity classes. 
In the last Section 4 we give further applications of the main theorem in Proposition~\ref{crMneN} that characterizes the matchings 
$M,N$ such that $cr({\cal T}(M,l))=ne({\cal T}(N,l))$ for every $l\ge 0$, in Corollary~\ref{camel} that concerns the statistic 
of pairs of separated edges, and in Theorem~\ref{modtwo} that enumerates the classes of mod $2$
crossing-similarity and mod $2$ nesting-similarity. We also give some concluding comments.

\section{The proof of Theorem~\ref{generalthm}}
For a set $X$ let $\M(X)$ be the set of all {\em finite multisets} with elements in $X$. By the sum
$$
X_1+X_2+\cdots+X_r=\sum_1^r X_i
$$
of the multisets $X_1,X_2,\dots,X_r\in\M(X)$ we mean the union of their groundsets with multiplicities of the elements added. 
Any function $f\colon X\to\M(Y)$ naturally extends to
$$
f:\;\M(X)\to\M(Y)\ \mbox{ by }\ f(U)=\sum_{x\in U}f(x) 
$$
where the summand $f(x)$ appears with the multiplicity of $x$ in $U$. Now if $Z\subset X$, we can understand the symbol $f(Z)$ 
in two ways---as the image of $f|Z$ or as the value of the extended $f$ on $Z$. Due to the convention on image we get in both cases 
the same result.

In this section $A$ shall denote an abelian group $(A,+)$ and $A^*$ will be the set of finite sequences over $A$. 
We shall work with functions 
from $A^*$ to $\M(A)$ or to $\M(A^*)$ which we will extend in the mentioned way, often without explicit notice, 
to functions defined on $\M(A^*)$. If $u=x_1x_2\dots x_t\in A^*$ and $y\in A$, by $x_1x_2\dots x_t+y$ we denote the sequence 
$(x_1+y)(x_2+y)\dots(x_t+y)$ obtained by adding $y$ to each term of $u$.

\begin{defi}\label{defini}
For $\alpha,\beta\in A$ and $i\in\N$ we define the mapping 
$R_{\alpha,\beta,i}:\;\bigcup_{l\ge i}A^l\to\bigcup_{l\ge i+2}A^l$ by
$$
R_{\alpha,\beta,i}(x_1x_2\dots x_l)=x_i(x_1x_2\dots x_i+x_i-x_1+\alpha)(x_ix_{i+1}\dots x_l+x_i-x_1+\beta)
$$
and the mapping $R_{\alpha,\beta}:\;A^*\to{\cal S}(A^*)$ by
$$
R_{\alpha,\beta}(x_1x_2\dots x_l)=\{R_{\alpha,\beta,i}(x_1x_2\dots x_l):\;1\le i\le l\}.
$$
So $R_{\alpha,\beta}(x_1x_2\dots x_l)$ is an $l$-element multiset of sequences which have length $l+2$. 
\end{defi}
 
Let $M\in{\cal M}(n)$ be a matching. The {\em gaps} of $M$ are the first gap before $[2n]$, the $2n-1$ 
gaps between the elements in $[2n]$, and the last $(2n+1)$-th gap after $[2n]$; $M$ has $2n+1$ gaps. 
For $\alpha,\beta\in A$ we assign to every matching $N\in{\cal M}(n)$, $n\in\N_0$, a sequence 
$seq_{\alpha,\beta}(N)\in A^*$ with length $2n+1$. If $n=0$, we set
$seq_{\alpha,\beta}(\emptyset)=0=0_A$. Let $n\ge 1$ and $(M,N)\in E({\cal T})$, $M\in{\cal M}(n-1)$, which means that $N$ is 
obtained from $M$ by adding a new first edge $e=\{1,x\}$ where $x$ is inserted in the $i$-th gap of $M$ for some $i$ lying  
between $1$ and $2n-1$. 
We set 
$$
seq_{\alpha,\beta}(N)=R_{\alpha,\beta,i}(seq_{\alpha,\beta}(M)).
$$
For example, if $M=\{\{1,3\},\{2,4\}\}$ then $seq_{\alpha,\beta}(M)=\alpha,2\alpha,3\alpha,2\alpha+\beta,\alpha+2\beta$.

For $u\in A^*$ we denote by $R^l_{\alpha,\beta}(u)=R_{\alpha,\beta}(R_{\alpha,\beta}(\dots(R_{\alpha,\beta}(u))\dots))$ the $l$-th 
iteration of the mapping $R_{\alpha,\beta}$ (which we extend to $\M(A^*)$). The next lemma is immediate from the definitions.

\begin{lema}\label{iterates}
For every $\alpha,\beta\in A$, $M\in{\cal M}$, and $l\in\N_0$ we have
$$
R^l_{\alpha,\beta}(seq_{\alpha,\beta}(M))=seq_{\alpha,\beta}({\cal T}(M,l)).
$$
\end{lema}

The next lemma relates the sequences $seq_{\alpha,\beta}(M)$ and the statistic $s_{\alpha,\beta}$ on ${\cal M}$.
\begin{lema}\label{firstterm}
For every $\alpha,\beta\in A$ and $N\in{\cal M}(n)$ the first term of the sequence $seq_{\alpha,\beta}(N)$ equals 
$s_{\alpha,\beta}(N)=cr(N)\alpha+ne(N)\beta$.
\end{lema}
\duk
For $n=0$ this holds. For $n\ge 1$ we proceed by induction on $n$.
Suppose that $(M,N)\in E({\cal T})$ and that $N$ arises by adding new first edge $\{1,x\}$ to $M$, where $x$ is inserted in 
the $i$-th gap. Let $seq_{\alpha,\beta}(M)=a_1a_2\dots a_{2n-1}$. 

We claim that in 
$$
a_j-a_1=u_j\alpha+v_j\beta
$$ 
the number $u_j$ counts the edges in $M$ covering the $j$-th gap and $v_j$ counts the edges in $M$ lying to the left of the $j$-th gap. 

Suppose that this claim holds. Then $cr(N)=cr(M)+u_i$ and $ne(N)=ne(M)+v_i$. Since $cr(M)\alpha+ne(M)\beta=a_1$ (by induction), 
the first term of $seq_{\alpha,\beta}(N)$ is $a_i=a_i-a_1+a_1=u_i\alpha+v_i\beta+cr(M)\alpha+ne(M)\beta=cr(N)\alpha+ne(N)\beta$, 
as we wanted to show.

It suffices to prove by induction on $n$ the claim. For $n=0$ it holds trivially. We assume that it holds for 
$seq_{\alpha,\beta}(M)$ and deduce it for $seq_{\alpha,\beta}(N)$; $M$, $N$, and $i$ are as before. Let 
$seq_{\alpha,\beta}(N)=b_1b_2\dots b_{2n+1}$. We first describe the changes in gaps 
caused by the addition of $\{1,x\}$ to $M$. A new first gap appears; it is of course covered by no edge and has no edge to its left. For 
$1\le j\le i$ the $j$-th gap turns in the $(j+1)$-th one; these gaps get covered by one more edge and have the same numbers 
of edges to their left as before. The $i$-th gap is split in two which creates a new gap, the $(i+2)$-th one; it is covered by 
as many edges as the $i$-th gap in $M$ but it has one more edge to its left. For 
$i+1\le j\le 2n-1$ the $j$-th gap turns in the $(j+2)$-th one; these gaps are covered by as may edges as before but they have 
one more edge to their left. 

By the definition of $R_{\alpha,\beta,i}$, $b_1=a_i$, 
$b_j=a_{j-1}+a_i-a_1+\alpha$ for $2\le j\le i+1$, and $b_j=a_{j-2}+a_i-a_1+\beta$ for $i+2\le j\le 2n+1$. Thus $b_1-b_1=0$, 
$b_j-b_1=a_{j-1}-a_1+\alpha=(u_{j-1}+1)\alpha+v_{j-1}\beta$ for $2\le j\le i+1$, and 
$b_j-b_1=a_{j-2}-a_1+\beta=u_{j-2}\alpha+(v_{j-2}+1)\beta$ for $i+2\le j\le 2n+1$. This agrees with the described changes in gaps 
and so the claim holds for $seq_{\alpha,\beta}(N)$.
\kduk

Let us denote by $f_0^0:\;A^*\to A$ the function taking the first term of a sequence and by $f^1_0:\;A^*\to \M(A)$ the function 
creating the multiset of all terms of a sequence. By the definitions and Lemmas~\ref{iterates} and \ref{firstterm}, if 
$seq_{\alpha,\beta}(M)=a_1a_2\dots a_{2n+1}$ then
$$
s_{\alpha,\beta}({\cal T}(M,1))=f_0^0(R_{\alpha,\beta}(seq_{\alpha,\beta}(M)))=\{a_1,a_2,\dots,a_{2n+1}\}=
f_0^1(seq_{\alpha,\beta}(M)).
$$
For the induction argument we will need more complicated functions besides $f_0^0$ and $f_0^1$. For an integer $r\ge 0$ and 
$\gamma\in A$ we define the function $f_{\gamma}^r\colon A^*\to\M(A)$ by
$$
f_{\gamma}^r(x_1x_2\dots x_l)=\{x_{a_1}+x_{a_2}+\cdots+x_{a_r}-(r-1)x_1+\gamma\colon
1\le a_1\le a_2\le\dots \le a_r\le l\}.
$$
So $f^0_0(x_1x_2\dots x_l)=\{x_1\}$ and $f^1_{\gamma}(x_1x_2\dots x_l)$ is the multiset 
$\{x_1+\gamma,x_2+\gamma,\dots, x_l+\gamma\}$.

\begin{lema}\label{main}
Let $X,Y\in\M(A^*)$ (possibly $X=Y$) be two multisets  such that 
$f_{\gamma}^r(X)=f_{\gamma}^r(Y)$ for every $r\ge 0$ and $\gamma\in A$. Then for every mapping $R_{\alpha,\beta}$ 
of Definition~\ref{defini} we have
\begin{enumerate}
\item $f_{\gamma}^r(R_{\alpha,\beta}(X))=f_{\gamma}^r(R_{\alpha,\beta}(Y))$\\
\item $f_{\gamma}^r(R_{\alpha,\beta}(X))=f_{\gamma}^r(R_{\beta,\alpha}(Y))$
\end{enumerate}
for every $r\ge 0$ and $\gamma\in A$. 
\end{lema}
\duk
We prove only the second identity with $R_{\alpha,\beta}$ and $R_{\beta,\alpha}$; the proof of the first one is similar and easier.
We proceed by induction on $r$. The case $r=0$ is clear since $f_{\gamma}^0(R_{\alpha,\beta}(X))=f_{\gamma}^1(X)$ 
for every $X\in\M(A^*)$ and $\gamma\in A$. We assume that $r\ge 1$ 
and that for every $s$, $0\le s<r$, and $\gamma\in A$ we have $f_{\gamma}^s(R_{\alpha,\beta}(X))=f_{\gamma}^s(R_{\beta,\alpha}(Y))$. We 
consider only the function $f_0^r$, the proof for general $\gamma$ is similar.

We split the multisets $U=f_0^r(R_{\alpha,\beta}(X))$ and $V=f_0^r(R_{\beta,\alpha}(Y))$, which arise by summation, in 
several contributions and show that after rearranging, the corresponding contributions to $U$ and $V$ are equal. 
$U$ is the multiset of elements $y_{a_1}+y_{a_2}+\cdots+y_{a_r}-(r-1)y_1$ where the sequence 
$y_1y_2\dots y_l$ runs through $R_{\alpha,\beta}(X)$ and the indices $a_i$ run through the 
$r$-tuples $1\le a_1\le a_2\le\dots\le a_r\le l$, and
similarly for $V$. The first contribution $C$ is defined by the condition $a_1=1$. $C$ contributes 
to $U$ the elements
$$
y_1+y_{a_2}+\cdots+y_{a_r}-(r-1)y_1=y_{a_2}+\cdots+y_{a_r}-(r-2)y_1, 
$$
where $y_1y_2\dots y_l$ runs through $R_{\alpha,\beta}(X)$
and the indices $a_i$ run through the $(r-1)$-tuples 
$1\le a_2\le a_3\le\dots\le a_r\le l$. Thus $C$ contributes $f_0^{r-1}(R_{\alpha,\beta}(X))$. To $V$ it 
contributes $f_0^{r-1}(R_{\beta,\alpha}(Y))$. Hence $C$ contributes equally to $U$ and $V$ because 
$f_0^{r-1}(R_{\alpha,\beta}(X))=f_0^{r-1}(R_{\beta,\alpha}(Y))$ by the inductive assumption. 

Each $v=y_1y_2\dots y_l\in R_{\alpha,\beta}(X)$ is in $R_{\alpha,\beta}(u)$ for some $u=x_1x_2\dots x_{l-2}\in X$ and 
(by the definition of $R_{\alpha,\beta}$) consists of 
three segments: it starts with a term $x_i$ of $u$, then it comes $x_1\dots x_i$ 
termwise incremented by $x_i-x_1+\alpha$, and the third segment of $v$ is $x_i\dots x_{l-2}$ termwise 
incremented by $x_i-x_1+\beta$; similarly for $v\in R_{\beta,\alpha}(Y)$. We split the rest of $U$ and $V$ 
(in which $a_1>1$, i.e., every $y_{a_i}$ lies in the second or in the third segment) in $r+1$ disjoint 
contributions $C_t$ according to the number $t$, $0\le t\le r$, of the $y_{a_i}$'s lying in the second segment. 
By the definition of $R_{\alpha,\beta}$, $C_t$ contributes to $U$ the elements 
\begin{eqnarray*}
&&x_{b_1}+\cdots+x_{b_r}+t(x_i-x_1+\alpha)+(r-t)(x_i-x_1+\beta)-(r-1)x_i\\
&=&x_{b_1}+\cdots+x_{b_r}+x_i-rx_1+t\alpha+(r-t)\beta
\end{eqnarray*}
where $u=x_1x_2\dots x_{l-2}$ runs through $X$, the indices $b_j$ run through the $r$-tuples satisfying
$1\le b_1\le\dots\le b_t\le i\le b_{t+1}\le\dots\le b_r\le l-2$, and $i$ runs through $1\le i\le l-2$. 
(The length $l-2$ depends on $u$.) Effectively the indices $b_j$ and $i$ run through all weakly increasing 
$(r+1)$-tuples of numbers from $[l-2]$. Thus $C_t$ contributes to $U$ the elements 
$f^{r+1}_{\gamma}(X)$ where $\gamma=t\alpha+(r-t)\beta$. By the definition of $R_{\beta,\alpha}$, $C_t$ contributes to 
$V$ the elements $f^{r+1}_{\gamma'}(Y)$ where $\gamma'=t\beta+(r-t)\alpha$. So $C_t$ contributes to $U$ and $V$ in general 
differently but (by the assumption on $X$ and $Y$) the contributions of $C_t$ to $U$ and $C_{r-t}$ 
to $V$ are equal. By symmetry, $\sum_0^r C_i$ contributes the same amount to $U$ and $V$. 
Since $U$ and $V$ are covered by equal and disjoint contributions $C$ and $\sum_0^r C_i$, we conclude that 
$U=V$, i.e., $f_0^r(R_{\alpha,\beta}(X))=f_0^r(R_{\beta,\alpha}(Y))$. 

The proof of 1 is similar and easier, because now $C_t$ contributes equally to 
$U=f_0^r(R_{\alpha,\beta}(X))$ and $V=f_0^r(R_{\alpha,\beta}(Y))$.
\kduk

Next we show that for the equality of all functions $f^r_{\gamma}$ on two one-element sets it in fact suffices that 
$f^0_0$ and $f^1_0$ are equal. We prove it in two lemmas. Let $g^r\colon A^*\to\M(A)$ be defined by
$$
g^r(x_1x_2\dots x_l)=\{x_{a_1}+x_{a_2}+\cdots+x_{a_r}\colon 1\le a_1\le a_2\le\dots \le a_r\le l\}.
$$

\begin{lema}\label{gr}
If $u,v\in A^*$ are such that $g^1(u)=g^1(v)$ then $g^r(u)=g^r(v)$ for all $r\ge 1$.
\end{lema}
\duk
Let $g^1(u)=g^1(v)$ and $r\in\N$. For $\bar{a}=(a_1,\dots,a_s)\in A^s$ we denote 
$(\bar{a})$ the multiset $\{a_1,\dots,a_s\}$ and if $\bar{n}=(n_1,\dots,n_s)\in\N^s$ then   
$\bar{n}\cdot\bar{a}=n_1a_1+\cdots+n_sa_s\in A$. For $s\in\N$, $X\in\M(A)$, and $u=x_1x_2\dots x_l\in A^*$ we 
denote 
$$
S(s,X,u)=\{\bar{x}=(x_{a_1},\dots,x_{a_s})\colon 1\le a_1<a_2<\dots<a_s\le l,(\bar{x})=X\}.
$$
For $r,s\in\N$ we denote
$$
N(r,s)=\{(n_1,\dots,n_s)\in\N^s\colon n_1+\cdots+n_s=r\}.
$$
Now we can rewrite $g^r(u)$ and $g^r(v)$ as 
\begin{eqnarray*}
g^r(u)&=&\{\bar{n}\cdot\bar{a}:\ s\in[r],X\in\M(A),\bar{n}\in N(r,s),\bar{a}\in S(s,X,u)\}\\
g^r(v)&=&\{\bar{n}\cdot\bar{a}:\ s\in[r],X\in\M(A),\bar{n}\in N(r,s),\bar{a}\in S(s,X,v)\}.
\end{eqnarray*}
We claim that (i) for every fixed $s\in[r]$ and $X\in\M(A)$ the multiset
$$
m(\bar{a})=\{\bar{n}\cdot\bar{a}:\ \bar{n}\in N(r,s)\}
$$
is the same for all $\bar{a}\in A^s$ with $(\bar{a})=X$ 
and that (ii) for every fixed $s\in[r]$ and $X\in\M(A)$ we have $|S(s,X,u)|=|S(s,X,v)|$. This will prove that 
$g^r(u)=g^r(v)$. 

To show (i), we take $\bar{a},\bar{b}\in A^s$ with $(\bar{a})=(\bar{b})=X$. Then $\bar{a}$ can be obtained from 
$\bar{b}$ by permuting coordinates: $\bar{a}=\pi(\bar{b})$ for some $\pi\in{\cal S}_s$, and
$\bar{n}\cdot\bar{b}=\pi(\bar{n})\cdot\bar{a}$. If $\bar{n}$ runs through $N(r,s)$, so does $\pi(\bar{n})$. Hence
$m(\bar{a})=m(\bar{b})$. To show (ii), we suppose that $X$ consists of the distinct elements $x_1,\dots,x_t$ with 
multiplicities $n_1,\dots,n_t$ where $n_1+\cdots+n_t=s$ (else $|S(s,X,u)|=|S(s,X,v)|=0$) and denote by $m_a(u)$ and 
$m_a(v)$ the numbers of occurrences of $a\in A$ in $u$ and $v$. Because $m_a(u)=m_a(v)$ for every $a\in A$, we have indeed 
$$
|S(s,X,u)|=\prod_{i=1}^t{m_{x_i}(u)\choose n_i}=\prod_{i=1}^t{m_{x_i}(v)\choose n_i}=|S(s,X,v)|.
$$
\kduk

\begin{lema}\label{just01}
If $X,Y\in\M(A^*)$ are one-element sets such that $f_0^0(X)=f_0^0(Y)$ and $f_0^1(X)=f_0^1(Y)$,
then $f_{\gamma}^r(X)=f_{\gamma}^r(Y)$ for every $r\ge 0$ and $\gamma\in A$.
\end{lema}
\duk 
We need to prove that if $u,v\in A^*$ are two sequences beginning with the same term and
having equal numbers of occurrences of each $a\in A$, then $f_{\gamma}^r(u)=f_{\gamma}^r(v)$ for every $r\ge 0$ 
and $\gamma\in A$. It suffices to consider functions $f^r_0$, the proof with general $\gamma$ is similar. 
Since $u$ and $v$ start with the same term, by the definition of $f^r_0$ it suffices to prove that $g^r(u)=g^r(v)$ 
for every $r\ge 1$. This is true by Lemma~\ref{gr}.
\kduk

\bigskip\noindent
{\bf Proof of Theorem~\ref{generalthm}.} We prove only 2, the proof of 1 is very similar and easier. Let 
$s_{\alpha,\beta}({\cal T}(M,l))=s_{\beta,\alpha}({\cal T}(N,l))$ for $l=0,1$. By Lemma~\ref{firstterm} and the following remark, 
this means that $f^0_0(seq_{\alpha,\beta}(M))=f^0_0(seq_{\beta,\alpha}(N))$ and 
$f^1_0(seq_{\alpha,\beta}(M))=f^1_0(seq_{\beta,\alpha}(N))$.
By Lemma~\ref{just01}, $f^r_{\gamma}(seq_{\alpha,\beta}(M))=f^r_{\gamma}(seq_{\beta,\alpha}(N))$
for every $r\in\N_0$ and $\gamma\in A$. By repeated application of 2 of Lemma~\ref{main} we get 
$$
f^r_{\gamma}(R^l_{\alpha,\beta}(seq_{\alpha,\beta}(M)))=f^r_{\gamma}(R^l_{\beta,\alpha}(seq_{\beta,\alpha}(N)))
$$
for every $l,r\in\N_0$ and $\gamma\in A$. In particular, 
$$
f^0_0(R^l_{\alpha,\beta}(seq_{\alpha,\beta}(M)))=f^0_0(R^l_{\beta,\alpha}(seq_{\beta,\alpha}(N))).
$$
But by Lemma~\ref{iterates} we have 
$$
R_{\alpha,\beta}^l(seq_{\alpha,\beta}(M))=seq_{\alpha,\beta}({\cal T}(M,l))\ \mbox{ and }\ 
R_{\beta,\alpha}^l(seq_{\beta,\alpha}(N))=seq_{\beta,\alpha}({\cal T}(N,l)).
$$  
Thus, by Lemma~\ref{firstterm}, 
$$
s_{\alpha,\beta}({\cal T}(M,l))=s_{\beta,\alpha}({\cal T}(N,l))
$$
for every $l\ge 0$, which we wanted to prove.
\kduk

We give a formulation of Theorem~\ref{generalthm} in terms of the sequences $seq_{\alpha,\beta}(M)$.

\begin{veta}\label{reform}
Let $M,N\in{\cal M}(n)$ be two (not necessarily distinct) matchings and $\alpha,\beta\in A$ be two elements of the abelian groups.
\begin{enumerate}
\item We have $s_{\alpha,\beta}({\cal T}(M,l))=s_{\alpha,\beta}({\cal T}(N,l))$ for all $l\ge 0$ iff 
$s_{\alpha,\beta}(M)=s_{\alpha,\beta}(N)$ and the sequences $seq_{\alpha,\beta}(M)$ and $seq_{\alpha,\beta}(N)$ are equal as 
multisets (when order is neglected).
\item We have $s_{\alpha,\beta}({\cal T}(M,l))=s_{\beta,\alpha}({\cal T}(N,l))$ for all $l\ge 0$ iff 
$s_{\alpha,\beta}(M)=s_{\beta,\alpha}(N)$ and the sequences $seq_{\alpha,\beta}(M)$ and $seq_{\beta,\alpha}(N)$ are equal as 
multisets.
\end{enumerate}

\end{veta}

\section{The numbers of similarity classes}

In this section we determine the cardinalities $|{\cal M}(n)/\!\sim_{cr}\!|$ and $|{\cal M}(n)/\!\sim_{ne}\!|$. Let $A=(\Z,+)$. 
For $M\in{\cal M}$ we define its {\em crossing sequence} $crs(M)$ by $crs(M)=seq_{1,0}(M)-a_1$, where $a_1$ is the first term of 
$seq_{1,0}(M)$, and its {\em nesting sequence} $nes(M)$ by $nes(M)=seq_{0,1}(M)-b_1$, where $b_1$ is the first term of 
$seq_{0,1}(M)$. Recall that (by the proof of Lemma~\ref{firstterm}) the $i$-th term of $crs(M)$ is the 
number of edges in $M$ covering the $i$-th gap and the $i$-th term of $nes(M)$ is the number of edges lying to the left of the 
$i$-th gap. For example, $M=\{\{1,4\},\{2,5\},\{3,6\}\}$ has $crs(M)=(0,1,2,3,2,1,0)$ and $nes(M)=(0,0,0,1,2,3)$. 
By Theorems~\ref{mainthm} and \ref{reform}, $M\sim_{cr}N\iff cr(M)=cr(N)\;\&\;f^1_0(crs(M))=f^1_0(crs(N))$, that is, 
$M$ and $N$ are crossing-similar iff they have the same numbers of crossings and their crossing sequences are equal as multisets; 
an analogous result holds for the nesting-similarity. 

Let $e=\{a,d\},f=\{b,c\}\in M$, $1\le a<b<c<d\le 2n$, be a nesting in $M\in{\cal M}(n)$. We define its {\em width} as 
$\min(b-a,d-c)$. We define the width of a crossing in the same way, only $\{a,d\}$ is replaced with $\{a,c\}$ and $\{b,c\}$ with 
$\{b,d\}$. Suppose the nesting $e,f$ has the minimum width among all nestings in $M$ and its width is realized by $b-a$. Switching 
the first vertices of the edges $e$ and $f$, we obtain another matching $N$. If the width of $e,f$ is realized by $d-c$, we switch 
the second vertices of $e$ and $f$. This transformation $M\leadsto N$ is called the {\em n-c transformation}. In the same way, 
by switching the first or the second vertices of the edges in a crossing with minimum width, we define the {\em c-n transformation}.

\begin{lema}\label{cnnctrans}
Let $M,N\in{\cal M}(n)$ where $N$ is obtained from $M$ by the n-c (c-n) transformation. Then $N$ has the same sets of first 
and second vertices of the edges as $M$ and $ne(N)=ne(M)-1,cr(N)=cr(M)+1$ ($ne(N)=ne(M)+1,cr(N)=cr(M)-1$). 
\end{lema}
\duk
The first claim about $N$ is obvious. Let $e=\{a,c\},f=\{b,d\}\in M$, $1\le a<b<c<d\le 2n$, be a crossing in $M$ with the minimum
width which is equal to $b-a$ (if it is equal to $d-c$, the argument is similar). The c-n transformation replaces $e$ by 
$e'=\{b,d\}$ and $f$ by $f'=\{a,d\}$. Because of the minimality of the width, every edge of $M$ that has one endpoint between $a$ and 
$b$ must have the other endpoint between $a$ and $b$ as well. It follows that $e'$ crosses the same edges distinct from $f$ as $e$ 
does and similarly for $f'$ and $f$. The edge $e'$ is covered by the same edges different from $f'$ as $e$ and similarly for
$f'$ and $f$. The edge $e'$ does not cover the edges lying between $a$ and $b$ which were covered by $e$ but these are now 
covered by $f'$ and were not covered by $f$. If we do not consider the pairs $e,f$ and $e',f'$, $M$ and $N$ have the same 
numbers of crossings and the same numbers of nestings. Since $e,f$ is a crossing and $e',f'$ is a nesting, in total $N$ has one less 
crossing and one more nesting than $M$. The argument for the n-c transformation is similar and is left to the reader.
\kduk

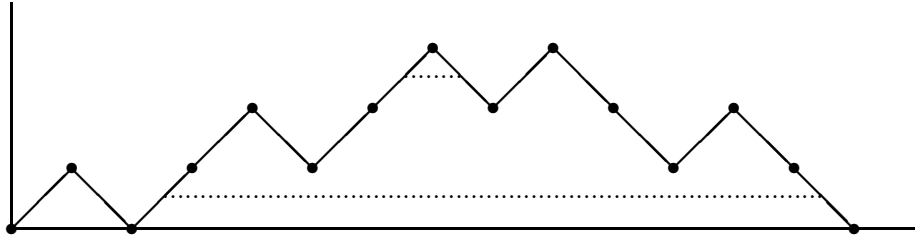
\begin{figure}
\unitlength1mm
\begin{picture}(60,40)(-5,0)
\put(0,0){\line(0,1){30}}\put(0,0){\line(1,0){120}}
\thicklines
\put(0,0){\line(1,1){8}}\put(8,8){\line(1,-1){8}}
\put(16,0){\line(1,1){8}}\put(24,8){\line(1,1){8}}\put(32,16){\line(1,-1){8}}
\put(40,8){\line(1,1){8}}\put(48,16){\line(1,1){8}}\put(56,24){\line(1,-1){8}}
\put(64,16){\line(1,1){8}}\put(72,24){\line(1,-1){8}}\put(80,16){\line(1,-1){8}}
\put(88,8){\line(1,1){8}}\put(96,16){\line(1,-1){8}}\put(104,8){\line(1,-1){8}}
\put(0,0){\circle*{1.5}}\put(8,8){\circle*{1.5}}\put(16,0){\circle*{1.5}}
\put(24,8){\circle*{1.5}}\put(32,16){\circle*{1.5}}\put(40,8){\circle*{1.5}}
\put(48,16){\circle*{1.5}}\put(56,24){\circle*{1.5}}\put(64,16){\circle*{1.5}}
\put(72,24){\circle*{1.5}}\put(80,16){\circle*{1.5}}\put(88,8){\circle*{1.5}}
\put(96,16){\circle*{1.5}}\put(104,8){\circle*{1.5}}\put(112,0){\circle*{1.5}}
\multiput(20,4)(1,0){88}{$.$} \multiput(52,20)(1,0){8}{$.$}
\end{picture}
\caption{Dyck path with semilength $7$ and two tunnels.}
\label{dyckovka}
\end{figure}

We use {\em Dyck paths} to encode $crs(M)$ and $nes(M)$. Recall that a Dyck path $D$ with semilength $n\in\N$ is a lattice 
path $D=(d_0,d_1,\dots,d_{2n})$, where $d_i\in\Z^2$, from $d_0=(0,0)$ to $d_{2n}=(2n,0)$ that makes $n$ up-steps 
$d_i-d_{i-1}=(1,1)$, $n$ down-steps $d_i-d_{i-1}=(1,-1)$, and never gets below the $x$ axis (so in fact $d_i\in\N_0^2$).
We denote the set of Dyck paths with semilength $n$ by ${\cal D}(n)$; $|{\cal D}(0)|=1$. We think of $D\in{\cal D}(n)$ also as 
a broken line in the 
plane that connects $(0,0)$ with $(2n,0)$ and consists of $2n$ straight segments $s_i=d_id_{i+1}$, see Figure~\ref{dyckovka}. 
A {\em tunnel} in $D$ is a horizontal segment $t$ that has altitude $n+\frac{1}{2}$ for some $n\in\N_0$, lies below $D$, and 
intersects $D$ only in its endpoints. Each $D\in{\cal D}(n)$ has exactly $n$ tunnels. Note that projections of two tunnels on 
the $x$ axis are either disjoint or they are in inclusion (as in the example in Figure~\ref{dyckovka}). If the latter happens, we 
say that the tunnel with larger projection {\em covers} the other tunnel. 

Deleting from $D\in{\cal D}(n)$, $n\ge 1$, the first up step and the first downstep at which $D$ visits again the $x$ axis,
we obtain, shifting appropriately the resulting two parts of $D$, a unique decomposition of $D$ in a pair of Dyck paths $E,F$, 
where $E\in{\cal D}(m)$ for $0\le m<n$ and $F\in{\cal D}(n-1-m)$. This {\em decomposition of Dyck paths} 
can be used for inductive proofs of their properties.
 
We associate with every Dyck path $D=(d_0,d_1,\dots,d_{2n})$ its {\em sequence of altitudes} 
$als(D)=(d_0^y,d_1^y,\dots,d_{2n}^y)\in\N_0^{2n+1}$, where $d_i=(d_i^x,d_i^y)$, and its {\em profile} 
$pr(D)=(a_1,a_2,\dots,a_m)\in\N^m$, where $m$ is the maximum term of $als(D)$ and $a_i$ is half of the number of segments $s_i$ 
of $D$ that lie in the horizontal strip $i-1\le y\le i$. It follows that $a_1+a_2+\cdots+a_m=n$ and $pr(D)$ is a composition of $n$.
It follows easily by induction on $m$ that for every composition $a=(a_1,a_2,\dots,a_m)$ of $n$ there is a $D\in{\cal D}(n)$ with
$pr(D)=a$. For example, the Dyck path in Figure~\ref{dyckovka} has $als(D)=(0,1,0,1,2,1,2,3,2,3,2,1,2,1,0)$ and 
$pr(D)=(2,3,2)$.

There is a natural surjective mapping $F:\;{\cal M}(n)\to{\cal D}(n)$ defined as follows. We take the diagram of $M\in{\cal M}(n)$ 
and travel the baseline $l$ from $-\infty$ to $\infty$. Simultaneously we construct, 
step by step, a lattice path $D$. We start $D$ at $(0,0)$ and when we encounter on $l$ the first (second) vertex of an edge, 
we make 
in $D$ an up-step (down-step). In the end we get a Dyck path $D\in{\cal D}(n)$ and set $F(M)=D$. Using the 
decomposition of Dyck paths and induction, it is easy to prove that $F$ is surjective. Clearly, the preimages $F^{-1}(D)$ 
consist exactly of the matchings sharing the same sets of first and second vertices. Another important property of $F$ is that 
for every $D\in{\cal D}(n)$ there is exactly one {\em noncrossing} (i.e., with $cr(M)=0$) $M\in F^{-1}(D)$, namely the $M$ 
whose edges correspond in the obvious way to the tunnels in $D$. This follows by the decomposition of Dyck paths.

\begin{lema}\label{crclasses}
Let $n\in\N$ and $F:\;{\cal M}(n)\to{\cal D}(n)$ be the above mapping.
\begin{enumerate} 
\item For every $M\in{\cal M}(n)$ we have $crs(M)=als(F(M))$. 
\item For every $M,N\in{\cal M}(n)$ we have $f^1_0(crs(M))=f^1_0(crs(N))$ iff $pr(F(M))=pr(F(N))$. 
\item For every composition $a=(a_1,a_2,\dots,a_m)$ of $n$ and every $i\in\N_0$, 
$0\le i\le\sum_{i=1}^m (i-1)a_i$, there is an $M\in{\cal M}(n)$ such that $pr(F(M))=a$ and $cr(M)=i$. There exist no 
$a$ and no $M$ such that $pr(F(M))=a$ and $cr(M)>\sum_{i=1}^m (i-1)a_i$.
\end{enumerate}
\end{lema}
\duk
1. This is clear from the definitions of $crs(M)$ and $als(D)$.

2. Using 1, we look at $f^1_0(als(D))$ where $D=F(M)$. Let $pr(D)=(a_1,a_2,\dots,a_m)$ and 
$r_i$ be the multiplicity of $i\in\N_0$ in $als(D)$. It is clear that
$r_0=a_1+1$ and $r_m=a_m$. We claim that for $0<i<m$ we have $r_i=a_i+a_{i+1}$. In the strip $i-1\le y\le i$ we have 
$v=2a_i$ segments $s_1,s_2,\dots,s_v$ of $D$ and in the strip $i\le y\le i+1$ we have $w=2a_{i+1}$ segments $t_1,t_2,\dots, t_w$.
The occurrences of $i$ in $als(D)$ are due to the upper endpoints of the $s_j$'s and due to the lower endpoints of the $t_j$'s. 
But for each $s_j$ its upper endpoint coincides with the upper endpoint of $s_{j-1}$ or with that of $s_{j+1}$ or with 
the lower endpoint of some $t_k$, and similarly for the lower endpoints of the $t_j$'s. So $i$ appears $(v+w)/2=a_i+a_{i+1}$ 
times. On the other hand, $a_i=r_{i-1}-r_{i-2}+\cdots+(-1)^{i}r_1+(-1)^{i+1}(r_0-1)$ for every $1\le i\le m$.
Therefore the $r_i$'s are completely determined by the composition $pr(D)$ and vice versa. 

3. Let a composition $a=(a_1,a_2,\dots,a_m)$ of $n$ be given. We take an arbitrary $D\in{\cal D}(n)$ with $pr(D)=a$. It follows by 
the decomposition of Dyck paths and induction that the sum
$$
S(a)=\sum_{i=1}^m (i-1)a_i
$$
counts the ordered pairs $t_1,t_2$ of distinct tunnels in $D$ where $t_1$ covers $t_2$. For the unique noncrossing 
$M\in F^{-1}(D)$ we have $ne(M)=S(a)$ because nestings in $M$ are in 1-1 correspondence with the pairs of tunnels, 
one of them covering the other. So $cr(M)=0,ne(M)=S(a),F(M)=D,pr(F(M))=a$. For any given $i\in\{0,1,\dots,S(a)\}$, 
using repeatedly the n-c transformation of Lemma~\ref{cnnctrans}, we transform $M$ into $N$ such that $cr(N)=i,ne(N)=S(a)-i,$ and 
$F(N)=F(M)=D$. Now suppose that there is an $M\in F^{-1}(D)$ with $cr(M)=c>S(a)$. Using the c-n transformation of 
Lemma~\ref{cnnctrans} we transform it into $N\in F^{-1}(D)$
with $cr(N)=0$ and $ne(N)=ne(M)+c>S(a)$. This contradicts the unicity of the noncrossing matching in $F^{-1}(D)$.
\kduk

\begin{veta}\label{cros_equi}
For $n\in\N$ the set ${\cal M}(n)/\!\sim_{cr}$ of crossing-similarity classes has
$$
2^{n-2}\bigg({n\choose 2}+2\bigg)
$$
elements.
\end{veta}
\duk
By the previous lemma, $|{\cal M}(n)/\!\sim_{cr}|$ equals
$$
\sum_a(1+a_2+2a_3+\cdots+(m-1)a_m)=2^{n-1}+\sum_a(a_2+2a_3+\cdots+(m-1)a_m)
$$
where we sum over all compositions $a_1+a_2+\cdots+a_m=n$, which are $2^{n-1}$ in number. The last sum is the coefficient 
of $x^n$ in the expansion of
$$
\bigg(\frac{d}{dy}\sum_{m\ge 0}\frac{x}{1-x}\cdot\frac{xy}{1-xy}\cdot\frac{xy^2}{1-xy^2}\cdot\dots\cdot\frac{xy^m}{1-xy^m}\bigg)
\bigg|_{y=1}.
$$
Differentiating the product in the summand by the Leibniz rule and using that 
$$
\bigg(\frac{d}{dy}\frac{xy^i}{1-xy^i}\bigg)\bigg|_{y=1}=\frac{ix}{(1-x)^2},
$$
we obtain that the expansion equals
$$
\frac{1}{1-x}\sum_{m\ge 0}{m+1\choose 2}\bigg(\frac{x}{1-x}\bigg)^{m+1}.
$$
Using the binomial expansion $(1-z)^{-r}=\sum_{n\ge 0}{r+n-1\choose n}z^n$, we simplify this to
$$
\frac{x^2}{(1-2x)^3}=\sum_{n\ge 0}{n+2\choose 2}2^nx^{n+2}
$$
and the result follows.
\kduk

\noindent
The values of $|{\cal M}(n)/\!\sim_{cr}\!|$ form the sequence $(1,3,10,32,96,276,\dots)$. Subtracting $2^{n-1}$, we get 
the sequence $(0,1,6,24,80,240,\dots)$ that counts crossing-similarity classes in ${\cal M}(n)$ for matchings with at least one 
crossing. This sequence is entry A001788 of \cite{sloa} and counts, for example, also $4$-cycles in the $(n+1)$-dimensional hypercube. 

The situation for nestings is simpler and the number of similarity classes is bigger because nesting sequences 
are nondecreasing and therefore $f^1_0(nes(M))=f^1_0(nes(N))$ iff $nes(M)=nes(N)$. By Theorems~\ref{mainthm} and 
\ref{reform}, $M\sim_{ne}N$ iff
$M$ and $N$ have the same numbers of nestings and the same nesting sequences. For $D\in{\cal D}(n)$ we define $ne(D)$ to be 
the number of ordered pairs $t_1,t_2$ of distinct tunnels in $D$ such that $t_1$ covers $t_2$. The {\em down sequence} 
$dos(D)$ of $D=(d_0,d_1,\dots,d_{2n})$ is $(v_0,v_1,\dots,v_{2n})$ where $v_i$ is the number of down-steps $d_j-d_{j-1}=(1,-1)$ 
for $1\le j\le i$. For example, for the Dyck path in Figure~\ref{dyckovka} we have $ne(D)=8$ and 
$dos(D)=(0,0,1,1,1,2,2,2,3,3,4,5,5,6,7)$.

\begin{lema}\label{dosseq}
Let $n\in\N$ and $F:\;{\cal M}(n)\to{\cal D}(n)$ be the mapping defined above.
\begin{enumerate} 
\item For every $M\in{\cal M}(n)$ we have $nes(M)=dos(F(M))$. There is a bijection between the sets 
$\{nes(M):\;M\in{\cal M}(n)\}$ and ${\cal D}(n)$.
\item For every Dyck path $D\in{\cal D}(n)$ and every $i\in\N_0$, $0\le i\le ne(D)$, there is an $M\in F^{-1}(D)$ such 
that $ne(M)=i$. There is no $M\in F^{-1}(D)$ with $ne(M)>ne(D)$.
\end{enumerate} 
\end{lema}
\duk
1. The first claim follows at once from the definitions. It is also clear that $dos(D)$ is uniquely determined by $D$ and vice 
versa. 

2. We know from the proof of 3 of Lemma~\ref{crclasses} that $ne(D)=ne(M)$ for the unique noncrossing $M\in F^{-1}(D)$. Now we 
argue as in the proof of 3 of Lemma~\ref{crclasses}.
\kduk

\begin{veta}\label{nest_equi}
For $n\in\N$ the set ${\cal M}(n)/\!\sim_{ne}$ of nesting-similarity classes has
$$
2\cdot 4^{n-1}-\frac{3n-1}{2n+2}{2n\choose n}
$$
elements.
\end{veta}
\duk
By the previous lemma,
$$
|{\cal M}(n)/\!\sim_{ne}\!|=\sum_{D\in{\cal D}(n)}(1+ ne(D))=|{\cal D}(n)|+\sum_{D\in{\cal D}(n)}ne(D).
$$
We claim that this number is equal to the coefficient of $x^n$ in the expansion of the expression
$$
C+x^2(2xC'+C)^2C
$$
where $C=C(x)=\sum_{n\ge 0}|{\cal D}(n)|x^n=1+x+2x^2+5x^3+\cdots$. It is well known that 
$C=(1-\sqrt{1-4x})/2x=\sum_{n\ge 0}\frac{1}{n+1}{2n\choose n}x^n$. Using the relations $xC^2-C+1=0$ and $2xCC'+C^2=C'$
we simplify the expression to 
$$
2C(x)+\frac{1/2}{1-4x}-\frac{3/2}{\sqrt{1-4x}}.
$$
Using the expansion of $C(x)$, geometric series, and $(1-4x)^{-1/2}=\sum_{n\ge 0}{2n\choose n}x^n$ we obtain the 
formula.

To establish the claim, recall that $\sum_{D\in{\cal D}(n)}ne(D)$ counts the triples $(D,t_1,t_2)$ where $D\in{\cal D}(n)$ and
$t_1,t_2$ are two distinct tunnels in $D$ such that $t_1$ covers $t_2$. Let the segments of $D$ supporting $t_i$ be 
$r_i$ (up-step) and $s_i$ (down-step). Let the lower endpoints of the segments $r_i$ ($s_i$) be $a_i$ ($b_i$) and their
upper endpoints be $a_i'$ ($b_i'$), $i=1,2$. The deletion of the interiors of the segments $r_1,s_1,r_2$, and $s_2$ 
splits $D$ in five lattice paths $L_1,\dots,L_5$ where 
$L_1$ starts at $(0,0)$ and ends in $a_1$, $L_2$ starts at $a_1'$ and ends at $a_2$, $L_3$ starts at $a_2'$ and ends at $b_2'$,
$L_4$ starts at $b_2$ and ends at $b_1'$, and $L_5$ starts at $b_1$ and ends at $(2n,0)$.
Each $L_i$ is nonempty but may be just a single lattice point. The concatenation $L_1L_5$, where $L_5$ is appropriately shifted 
so that $a_1$ and $b_1$ are identified in one distinguished point, is a Dyck path and similarly for $L_2L_4$ with $a_2$ and $b_2$
identified and distinguished. $L_3$ is a Dyck path by itself (after an appropriate shift). We see that the triples 
$(D,t_1,t_2)$ in question are in a 1-1 correspondence with the triples $(E_1,E_2,E_3)$ where $E_i\in{\cal D}(n_i)$, 
$n_i\in\N_0$, are such that $n_1+n_2+n_3=n-2$, and moreover $E_1$ and $E_2$ have one distinguished lattice point 
(out of $2n_1+1$, respectively $2n_2+1$, points). It follows that the number of the triples $(E_1,E_2,E_3)$ is the coefficient of 
$x^{n-2}$ in $(2xC'+C)^2C$, which
proves the claim. 
\kduk

\noindent
The values of $|{\cal M}(n)/\!\sim_{ne}\!|$ form the sequence $(1,3,12,51,218,926,\dots)$. Subtracting the Catalan numbers 
$C_n=|{\cal D}(n)|$, we get the sequence $(0,1,7,37,176,794,\dots)$ that counts nesting-similarity classes in ${\cal M}(n)$ 
for matchings with at least 
one nesting. This sequence is entry A006419 of \cite{sloa} and appears in Welsh and Lehman \cite[Table VIb]{wels_lehm} in 
enumeration of planar maps. We summarize this identity in the next proposition.

\begin{figure}
\unitlength1mm
\begin{picture}(60,50)(-18,-4)
\put(0,5){\oval(22,10)}\put(8,5){\circle{6}}\put(-11,5){\circle*{1.5}}\put(11,5){\circle*{2}}
\put(0,20){\oval(22,10)}\put(-8,20){\circle{6}}\put(-11,20){\circle*{1.5}}\put(11,20){\circle*{2}}
\put(38,5){\oval(22,10)}\put(52,5){\circle{6}}\put(27,5){\circle*{1.5}}\put(49,5){\circle*{2}}
\put(38,20){\oval(22,10)}\put(24,20){\circle{6}}\put(27,20){\circle*{1.5}}\put(49,20){\circle*{2}}
\put(76,5){\oval(22,10)}\put(84,5){\circle{6}}\put(81,5){\circle*{1.5}}\put(87,5){\circle*{2}}
\put(76,20){\oval(22,10)}\put(68,20){\circle{6}}\put(65,20){\circle*{1.5}}\put(71,20){\circle*{2}}
\put(102,12.5){\oval(10,25)}\put(102,0){\line(0,1){25}}\put(102,25){\circle*{2}}\put(102,0){\circle*{1.5}}
\end{picture}
\caption{All seven rooted and edge $2$-connected planar maps with two vertices and three faces.}
\label{entri}
\end{figure}
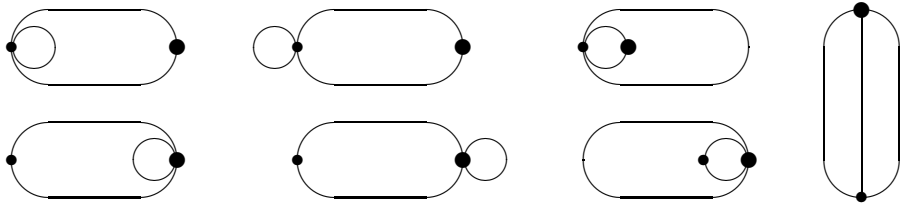

\begin{prop}
For $n=1,2,\dots$ the formula 
$$
2\cdot 4^{n-1}-\frac{3n+1}{2n+2}{2n\choose n}
$$
counts the following objects.
\begin{enumerate}
\item The triples $(D,t_1,t_2)$ where $D$ is a Dyck path with semilength $n$ and $t_1,t_2$ are two distinct tunnels in $D$ such that 
$t_1$ covers $t_2$.
\item The nesting-similarity classes in $\{M\in{\cal M}(n):\;ne(M)>0\}/\sim_{ne}$.
\item The vertex-rooted planar maps with two vertices and $n$ faces, which are edge $2$-connected and may have loops and multiple edges.
See Figure~\ref{entri} for the case $n=3$.
\end{enumerate}
\end{prop}
\duk
1 and 2 follow from the proof of Theorem~\ref{nest_equi} and 3 follows by checking the formulas in \cite{wels_lehm}. Alternatively, it 
is not too hard to establish bijection between the triples in 1 and the maps in 3. 
\kduk

 The present author
proved in \cite[Theorem 3.1]{klaz97} that the number of the triples $(T,v_1,v_2)$, where $T$ is a rooted plane tree with $n$ 
vertices and $v_1,v_2$ are two (not necessarily distinct) vertices of $T$ such that $v_1$ lies on the path joining the root 
of $T$ and $v_2$, equals
$$
\frac{4^{n-1}+{2n-2\choose n-1}}{2}.
$$ 
It is straightforward to relate Dyck paths and rooted plane trees and to derive the formula of Theorem~\ref{nest_equi} from this one.

\section{Further applications and concluding remarks}

Corollary~\ref{crequalne} presents two matchings $M$ and $N$ such that the distribution of $cr$ 
on the levels of ${\cal T}(M)$ equals the distribution of $ne$ on the levels of ${\cal T}(N)$. We show that there are no other 
substantially different examples.

\begin{prop}\label{crMneN}
Let $M,N\in{\cal M}(n)$ be two matchings. We have $cr({\cal T}(M,l))=ne({\cal T}(N,l))$ for every $l\ge 0$ if and only if 
$M=M_n=\{\{1,2n\},\{2,2n-1\},\dots,\{n,n+1\}\}$ and $N=N_n=\{\{1,2\},\{3,4\},\dots,\{2n-1,2n\}\}$.
\end{prop}
\duk
The if part is clear by Theorem~\ref{mainthm}: $cr(M_n)=ne(N_n)=0$ and 
$$
cr({\cal T}(M_n,1))=ne({\cal T}(N_n,1))=\{0,0,1,1,2,2,\dots,n-1,n-1,n\}
$$ 
because 
\begin{eqnarray*}
crs(M_n)&=&(0,1,2,\dots,n-1,n,n-1,\dots,2,1,0)\\
nes(N_n)&=&(0,0,1,1,2,2,\dots,n-1,n-1,n).
\end{eqnarray*}

To show the only if part, we prove that the only matchings $M,N\in{\cal M}(n)$ satisfying 
$cr(M)=ne(N)$ and $f^1_0(crs(M))=f^1_0(nes(N))$ are $M_n$ and $N_n$. Since for every $N\in{\cal M}(n)$ the sequence 
$nes(N)$ ends with $n$, we must have $n$ in $crs(M)$ which means that the middle gap of $M$ must be covered by all edges. Thus all 
first vertices of the edges in $M$ must precede all second vertices and $crs(M)=(0,1,2,\dots,n-1,n,n-1,\dots,2,1,0)$. Thus 
$f_0^1(nes(N))=\{0,0,1,1,2,2,\dots,n-1,n-1,n\}$ which forces $N=N_n$. Thus $cr(M)=ne(N)=ne(N_n)=0$ which forces $M=M_n$.
\kduk

\noindent
Therefore we have no other examples of equidistribution of $cr$ and $ne$ on the levels of ${\cal T}(M)$ than $M=\emptyset$ 
and $M=\{\{1,2\}\}$ because $M_n=N_n$ only for $n=0,1$. We call the matchings $M\in{\cal M}(n)$ 
encountered in the proof in which all edges cover the middle gap, equivalently which have 
$f_0^1(crs(M))=\{0,0,1,1,2,2,\dots,n-1,n-1,n\}$, {\em permutational matchings}; they are in 1-1 correspondence with the 
permutations of $[n]$ and are $n!$ in number. 

Because $|{\cal M}(n)|=(2n-1)!!=n^n(2/\mathrm{e}+o(1))^n$ and the numbers of crossing-similarity and nesting-similarity classes 
are only exponential, we have very many examples as in Corollary~\ref{crequalcr} when $cr$ (or $ne$) has equal distributions 
on the levels of ${\cal T}(M)$ and ${\cal T}(N)$ for $M\ne N$. The next corollary follows from the asymptotics of the numbers 
of similarity classes given in Theorems~\ref{cros_equi} and \ref{nest_equi}.

\begin{dusl}
Every set of matchings $X\subset{\cal M}(n)$ contains $|X|/(2+o(1))^n$ mutually crossing-similar matchings and 
$|X|/(4+o(1))^n$ mutually nesting-similar matchings.
\end{dusl}

\noindent
An explicit example of a big similarity class is provided by permutational matchings in ${\cal M}(n)$. They all share
the same crossing sequence $(0,1,2,\dots,n-1,n,n-1,\dots,2,1,0)$ and the same nesting sequence $(0,0,\dots,0,1,2,\dots,n-1,n)$. 
Hence at least
$$
\frac{n!}{{n\choose 2}+1}=n^n(1/\mathrm{e}+o(1))^n
$$
of them are mutually crossing-similar and at least so many of them are mutually nesting-similar.

Crossing and nesting correspond to two of three matchings in ${\cal M}(2)$ and the third remaining matching is 
$\{\{1,2\},\{3,4\}\}$. If two edges of $M\in{\cal M}$ form this matching, we say that they 
form a {\em camel}. We denote the number of camels in $M$ by $ca(M)$. This statistic behaves on the levels of the subtrees of 
${\cal T}$ in the same way as $cr$ and $ne$ do.

\begin{dusl}\label{camel}
Let $M,N\in{\cal M}(n)$ be two matchings such that $ca$ has the same distribution on the first two levels of the subtrees 
${\cal T}(M)$ and ${\cal T}(N)$. Then $ca$ has the same distribution on all levels.
\end{dusl}
\duk
For $M\in{\cal M}(n)$ we have $ca(M)={n\choose 2}-(cr(M)+ne(M))$. Thus this result follows by 1 of Theorem~\ref{generalthm} if we set
$A=(\Z,+)$ and $\alpha=\beta=1$.
\kduk

\noindent
Note that while the number of $M\in{\cal M}(n)$ with $cr(M)=0$ (or with $ne(M)=0$) is the Catalan number 
$\frac{1}{n+1}{2n\choose n}$, the number of $M\in{\cal M}(n)$ with $ca(M)=0$ is much bigger, namely $n!$ (these are exactly 
permutational matchings).

It is possible to investigate the general similarity relation $\sim_{A,\alpha,\beta}$ on ${\cal M}(n)$ defined, 
for an abelian group $A=(A,+)$ 
and two its elements $\alpha,\beta\in A$, by $M\sim_{A,\alpha,\beta}N$ iff 
$s_{\alpha,\beta}({\cal T}(M,l))=s_{\alpha,\beta}({\cal T}(N,l))$ for every $l\ge 0$. We consider here only the case $A=(\Z_2,+)$
and define the statistics $cr_2(M),ne_2(M)\in\{0,1\}$ as parity of the numbers $cr(M),ne(M)$. We define the
sequences $crs_2(M)$ and $nes_2(M)$ of $M$ by reducing $crs(M)$ and $nes(M)$ modulo $2$. For two matchings 
$M,N\in{\cal M}(n)$ we define $M\sim_{cr,2}N$ iff $cr_2({\cal T}(M,l))=cr_2({\cal T}(M,l))$ for every $l\ge 0$, and similarly 
for $M\sim_{ne,2}N$. By Theorems~\ref{generalthm} and \ref{reform}, $M\sim_{cr,2}N$ iff $cr_2(M)=cr_2(N)$ and 
$crs_2(M)$ and $crs_2(N)$ are equal as multisets after forgetting the order of terms, and similarly for $\sim_{ne,2}$. 
(Now $nes_2(M)$ is not nondecreasing and we may have $f^1_0(nes_2(M))=f^1_0(nes_2(N))$ for $nes_2(M)\ne nes_2(N)$.) We determine 
the numbers of equivalence classes for $\sim_{cr,2}$ and $\sim_{ne,2}$.

\begin{veta}\label{modtwo}
We have $|{\cal M}(1)/\!\sim_{cr,2}\!|=1$ and $|{\cal M}(n)/\!\sim_{cr,2}\!|=2$ for $n\ge 2$. The two classes of mod $2$ 
crossing-similarity  have $((2n-1)!!+1)/2$ and $((2n-1)!!-1)/2$ elements. We have $|{\cal M}(1)/\!\sim_{ne,2}\!|=1$, 
$|{\cal M}(2)/\!\sim_{ne,2}\!|=3$, and $|{\cal M}(n)/\!\sim_{ne,2}\!|=2n$ for $n\ge 3$.
\end{veta}
\duk
It follows from the definition of $crs(M)$ that $crs_2(M)=(0,1,0,1,0,\dots,1,0)$ for every matching $M$. 
Thus the classes of mod $2$ crossing-similarity are determined only by $cr_2(M)$ and, for $n\ge 2$, we have two of them. 
The fact that 
$$
|\{M\in{\cal M}(n):\;cr_2(M)=0\}|-|\{M\in{\cal M}(n):\;cr_2(M)=1\}|=1
$$
for every $n\ge 1$ was proved by Riordan \cite{rior} by generating functions; a simple proof by involution was given 
by Klazar \cite{klaz}.

To handle nestings modulo $2$, recall that $nes(M)=dos(D)$ where $D=F(M)$ and that nesting sequences of the matchings
$M\in{\cal M}(n)$ are in bijection with the Dyck paths $D\in{\cal D}(n)$ (Lemma~\ref{dosseq}). We claim that the $n$ Dyck paths
$$
D_1=udu^{n-1}d^{n-1},D_2=u^2du^{n-2}d^{n-1},\dots,D_{n-1}=u^{n-1}dud^{n-1},D_n=u^nd^n
$$
($u$ is the up-step and $d$ is the down-step) realize all possible numbers of $1$'s and $0$'s in the sequences
$\{dos_2(D):\;D\in{\cal D}(n)\}$ and hence in the sequences $\{nes_2(M):\;M\in{\cal M}(n)\}$. The number of $1$'s ($0$'s) in 
$dos_2(D_i)$, $i=1,2,\dots,n$, is $n+\lceil n/2\rceil-i$ ($1+i+\lfloor n/2\rfloor$). It suffices to show that 
no $dos_2(D)$ has fewer than $\lceil n/2\rceil$ $1$'s and fewer than 
$2+\lfloor n/2\rfloor$ $0$'s. In every $D$ each of the $n$ down-steps contributes to $dos_2(D)$ exactly one $1$ 
(by one of its endpoints) and each of these $1$'s may belong to at most two downsteps. So we must have at least 
$\lceil n/2\rceil$ $1$'s. The argument for $0$'s is similar, but now the $0$ contributed by the first down-step is never shared 
(with the next down-steps) and there is one more $0$ contributed by the first up-step. So we have at least 
$1+1+\lfloor n/2\rfloor$ $0$'s. Thus, for every $n\ge 1$, $|\{f^1_0(nes_2(M)):\;M\in{\cal M}(n)\}|=n$. If $n\ge 3$, for each 
$D_i$ there are $M,M'\in F^{-1}(D_i)$ with $ne(M')=ne(M)-1$ (we take for $M$ the noncrossing matching in $F^{-1}(D_i)$, it has 
least one nesting, and apply the n-c transformation). Thus, for $n\ge 3$, there are $2n$ classes of mod $2$ nesting-similarity. 
The cases $n=1,2$ are easy to treat separately.
\kduk

\bigskip\noindent
{\bf Concluding remarks.} Recently, an interesting result for crossings and nestings of higher order was obtained by 
Chen et al. in \cite{chen_etal} where it is proved that for every $k,l,n\in\N$ the number of matchings in ${\cal M}(n)$ with 
no $k$-crossing and no $l$-nesting is the same as the number of matchings with no $k$-nesting and no $l$-crossing (a similar result 
is in \cite{chen_etal} obtained for set partitions); here $k$-crossing is a $k$-tuple of pairwise crossing edges and 
similarly for $k$-nesting. Another generalization of crossings and nestings is investigated by Jel\'\i nek \cite{jeli} 
who is interested in numbers of matchings $M\in{\cal M}(n)$ such that $M$ does not contain a fixed permutational matching 
$N\in{\cal M}(3)$ as an ordered submatching. 

It may be interesting to try to extend results and methods of the present article to crossings and nestings of higher order. Another
research direction may be to apply our method to other structures besides matchings. Finally, one may try to go to higher levels 
of the description of the enumerative complexity of crossings and nestings --- denoting 
$G:\;{\cal M}\to{\cal M}/\sim_{cr}$ the mapping sending $M$ to its equivalence class, when is it the case that 
$G({\cal T}(M,l))=G({\cal T}(N,l))$ for every $l\ge 0$; and similarly for $\sim_{ne}$. 

\bigskip\noindent
{\bf Acknowledgments.} I am grateful to Marc Noy for his hospitality during my two visits in UPC Barcelona in 2004 and for 
many stimulating discussions.

\end{document}